\newmathalphabet*{\E}{eus}{m}{n}

\renewcommand{\P}{{\Bbb P}}
\newcommand{\A}{{\Bbb A}}
\newcommand{\F}{{\Bbb F}}
\newcommand{\G}{{\Bbb G}}
\renewcommand{\O}{{\E O}}
\newcommand{\Z}{{\Bbb Z}}
\newcommand{\N}{{\Bbb N}}

\newcommand{\<}{\langle}
\renewcommand{\>}{\rangle}

\newcommand{\red}{\operatorname{{red}}}
\newcommand{\kar}{\operatorname{char}}
\newcommand{\Spec}{\operatorname{Spec}}

\newcommand{\Lie}{\operatorname{Lie}}
\newcommand{\Dist}{\operatorname{Dist}}
\newcommand{\mod}{\operatorname{mod}}

\newcommand{\Mor}{\operatorname{Mor}}

\renewcommand{\H}{\operatorname{H}}

\newcommand{\Pic}{\operatorname{Pic}}

\newcommand{\SL}{\operatorname{SL}}

\newtheorem{example}{Example}
\newtheorem{thm}{Theorem}
\newtheorem{lem}{Lemma}
\newtheorem{prop}{Proposition}

\newenvironment{pf}{{\bf Proof}}{$\Box$}

\newcommand{\binom}[2]{{\left(\begin{array}{c} #1\\#2 \end{array}\right)}}

\documentstyle[amscd,amssymb,12pt]{amsart}

\title{Embeddings of Homogeneous Spaces in Prime Characteristics}

\author{Niels Lauritzen}
\address{Matematisk Institut\\
Aarhus Universitet\\
Ny Munkegade\\
DK-8000 \AA rhus C\\
Denmark}
\email{niels@@mi.aau.dk}



\begin{document}

\maketitle

Let $X$ be a projective algebraic variety over an algebraically
closed field $k$ admitting a homogeneous action of a semisimple
linear algebraic group $G$. Then $X$ can be canonically identified
with the homogeneous space $G/G_x$, where $x$ is a closed point in $X$ and
$G_x$ the stabilizer group scheme of $x$. A group scheme over
a field of characteristic $0$ is reduced so in this case, $X$ is
isomorphic to a generalized flag variety $G/P$, where
$P$ is a parabolic subgroup. In
\cite{HaLa:varunsep}\cite{Lau:Split}\cite{Lau:Linebundles}\cite{Lau:Euler} the
geometry of $X$ in prime characteristic has been studied and
it has been shown that a lot of strange phenomena occur
when $G_x$ is non-reduced. The simplest example of a projective homogeneous
$G$-space (for $G=\SL_3(k)$, $\kar\, k=p>0$) not isomorphic
to a generalized flag variety is the divisor
$x_0\, y_0^p + x_1\, y_1^p + x_2\, y_2^p=0$ in $\P^2\times \P^2$. Since
projective homogeneous spaces with non-reduced stabilizers are quite
algebraic by construction, we give in \S \ref{incvar} of this paper
a simple geometric approach for their construction, involving
only scheme theoretic images under partial Frobenius morphisms.
We choose to do this focusing on the ``unseparated incidence variety''. In
this case the geometric approach completely determines the cohomology
of effective line bundles.
The reader unfamiliar with the general concept of
projective homogeneous spaces in prime characteristic might
find this section useful.

%In \cite{La2} a counterexample to
%Kodaira's vanishing theorem was constructed using a carefully
%chosen non-reduced
%parabolic subgroup scheme of $\SO_{10}(k)$.

The main topic of this paper is the study of
embeddings of homogeneous projective spaces. Let $X$ be
a projective homogeneous $G$-space. In \S \ref{gaction} we show that
$X$ can be realized as the $G$-orbit
of the $B$-stable line in $\P(L(\lambda))$, where
$L(\lambda)$ denotes the simple $G$-representation of a certain highest
weight $\lambda$. This approach leads in \S \ref{ex} to examples of
some strange embeddings of homogeneous spaces in characteristic $2$ -
one lying on the boundary of Hartshorne's conjecture on complete
intersections (a socalled Hartshorne variety in the terms of
\cite{LazVen}).

Recall that an ample line bundle $L$ on an algebraic variety $X$ is called
normally generated \cite{Mum:Quadrequa} if the multiplication map
$$
\H^0(X, L)^{\otimes n}\rightarrow \H^0(X, L^{\otimes n})
$$
is surjective for all $n\geq 1$. For generalized flag varieties
the very ampleness of ample line bundles follows from normal
generation of ample line bundles \cite{And:Frobcoh}\cite{MeRam:Schubvar}
(In \cite{Mum:Quadrequa} one can
find a nice and short proof of the fact that an ample normally
generated line bundle $L$ is very ample).
In the
general setting of projective homogeneous spaces, normal
generation of ample line bundles is an open question.

In \S \ref{linebundles} we compute the line bundles on projective
homogeneous spaces and in
\S \ref{ample} we use a simple ``diagonal'' construction
to prove that
ample line bundles on projective homogeneous $G$-spaces are
very ample.
In view of \cite{Lau:Euler} this
proves the existence of a counterexample to Kodaira type vanishing in prime
characteristic with a very ample line bundle answering Raynaud's question
(\cite{Ray:Contreexample}, Remarques et questions 3).

%The embedding $X\hookrightarrow \P(L(\lambda))$ is
%rarely linearly normal.

%The main idea of the proof
%can be read off from the example in \S \ref{incvar}, where we give
%explicit non-trivial examples of projective homogeneous spaces, which
%is not (isomorphic to) a generalized flag variety. This example is
%already implicit in \cite{HaLa},
% but the underlying simple geometry,
%which was somewhat obscured by roots and weights,
%was pointed out to me by Kollar and Mehta.

%There is also a representation
%theoretic proof using embeddings into projective spaces on Weyl modules
%and the algebra of differential operators in \cite{Ja}.
%Let $L(\lambda)$
%be the simple $G$-module of highest weight $\lambda$.
%In \S
%\ref{Llambda} we compute the stabilizer $G_x$ of the $B$-stable line
%$x\in \P(L(\lambda))$ for the natural action of $G$. This leads
%to a result about the weight structure of $L(\lambda))$.

%When $X$ is a generalized
%flag variety normal generation of an ample line bundle follows
%from the fact that the diagonal $\Delta_X\subseteq X\times X$ is
%compatibly Frobenius split \cite{Ram}. There is also
%a representation theoretic proof of normal generation using properties
%of the Steinberg weight\cite{A}. All of these proofs break down
%in the general setting of projective homogeneous spaces, where
%normal generation of ample line bundles is an open question.
The paper is organized as follows

\tableofcontents

%vskip1truecm

It is a pleasure to acknowledge inspiring discussions with
and e-mails from H.~H.~Andersen, W.~J.~Haboush, J.~Koll\'ar,
V.~B.~Mehta and V.~Srinivas. I am indebted to
Koll\'ar and Mehta for pointing out the usefulness of the
simple geometric viewpoint in \S \ref{incvar} and to
Haboush for telling me about the beautiful little book
\cite{LazVen} containing Zak's classification of Severi
varieties.

Part of this work was done while the author
was visiting Max-Planck-Institut f\"ur Mathematik in Bonn. I thank the
Institute for its hospitality and an inspiring environment.

\section{Preliminaries}

\label{charassumption}

Let $k$ be an algebraically closed field of prime
characteristic $p$. In the following, an algebraic variety $X$ is
assumed to be an
algebraic variety over $k$ and a morphism to be a
morphism of $k$-varieties. We let $X(A)=\Mor_k(\Spec\, A, X)$ denote
the set of $A$-points of $X$, where $A$ is a $k$-algebra.
Let $G$ be a simply connected and
semisimple algebraic group. We will assume that
$p>3$ if $G$ has a component of
type $G_2$ and $p>2$ if $G$ has a component of type $B_n$, $C_n$
or $F_4$.

\subsection{$G$-spaces}

A $G$-space is an algebraic variety $X$ endowed with
a morphism $G\times X\rightarrow X$ inducing
an action of  $G(A)$ on $X(A)$ for all $k$-algebras $A$.
A $G$-space $X$ is called homogeneous if the action
$G(k)\times X(k)\rightarrow X(k)$ on $k$-points is
transitive. A $k$-point $x\in X(k)$ gives a natural
morphism $G\rightarrow X$. The fiber product $G_x=G\times_X\Spec(k)$
is easily seen to be a closed subgroup scheme of $G$. It is called
the stabilizer group scheme of $x$.

\subsection{The Frobenius subcover}

\label{fsubcover}

An algebraic variety $X$ gives rise to a new algebraic
variety $X^{(n)}$ with the same underlying point space
as $X$, but where the $k$-multiplication
is twisted via the ringhomomorphism: $a\mapsto \sqrt[p^n]{a}$. The
$n$-th order Frobenius homomorphism induces a natural morphism
$F^{n}_X: X\rightarrow X^{(n)}$. As $X$ is reduced, $\O_{X^{(n)}}$
can be identified with the $k$-subalgebra of $p^n$-th powers of
regular functions on $X$. We call $X^{(n)}$ the $n$-th Frobenius
subcover of $X$. Recall that $X$ is said to be defined over
$\F_p$ if there exists an $\F_p$-variety $X'$, such that
$X\cong X'\times_{\Spec\, \F_p} \Spec\, k$. If $X$ is defined
over $\F_p$, then $X$ is isomorphic to $X^{(n)}$ (the isomorphism
is given locally by $f\otimes a\mapsto f\otimes a^{p^n}$, where
$a\in k$).

\subsection{The Frobenius kernel}

\label{fkernel}

Now $G^{(n)}$ is an algebraic group and $F^n_G:G\rightarrow G^{(n)}$
is a homomorphism of algebraic groups. The kernel of $F^n_G$
is called the $n$-th Frobenius kernel of $G$ and denoted $G_n$.
Let $X$ be a homogeneous $G$-space and $x$ a closed point of $X$. It
is easy to see that $X^{(n)}$ is a $G$-homogeneous space through
the homomorphism $F^n_G$. If $H=G_x$ then $G_{F^n_G(x)}=G_n\, H$.

%(some more comments on this, please!).

\subsection{The diagonal action}

\label{product}

Now let $X$ and $Y$ be homogeneous $G$-spaces with distinguished
closed points $x$ and $y$ and let $H=G_x$ and $K=G_y$. The product
$X\times Y$ becomes a $G$-space through the diagonal action and
$G_{(x,y)}=H\cap K$.

\section{The unseparated incidence variety}

\label{incvar}

In this section we give a quite explicit geometric description (in
\ref{unsepdivisor}) of
certain projective homogeneous spaces for $\SL_n$ occuring
as divisors in $\P^n\times \P^n$.

\subsection{The incidence variety}

Let $n>1$ and $G=\SL_{n+1}(k)$. The natural action of $G$ on
$V=k^{n+1}$ makes $\P(V)$ and $\P(V^*)$ into homogeneous
spaces for $G$. We fix points $x_1\in \P(V)$ and $x_2\in \P(V^*)$,
such that $G_{x_1}=P_1$ and $G_{x_2}=P_2$
are appropriate parabolic subgroups containing the subgroup
of upper triangular matrices $B$ in $G$.
The orbit $Y$ of $(x_1, x_2)$ in $\P(V)\times \P(V^*)$ is a projective
homogeneous space for $G$ isomorphic to $G/P$, where $P=P_1\cap P_2$.
Notice that the points of $Y$ are just pairs of incident
lines and hyperplanes and that $Y=Z(s)$, where $s$ is
the section $x_0 y_0 +\dots + x_n y_n$ of $\O(1)\times\O(1)$.

\subsection{The unseparated incidence variety}

\label{unsepdivisor}

Let $X$ be the $G$-orbit of
of $(x_1, F^r(x_2))$ in $\P(V)\times \P(V^*)^{(r)}$.
By \ref{fkernel} and \ref{product},
$$
X\cong G/\tilde{P}
$$
where for any $k$-algebra $A$
$$
\tilde{P}(A)=\{\left(
\begin{array}{lllll}
*        & *        & *            & \dots        & *       \\
0        & *        & *            & \dots        & *       \\
0        & *        & *            & \dots        & *       \\
\vdots   & \vdots   & \vdots       & \vdots       & \vdots  \\
0        & *        & *            & \dots        & *       \\
0        & a        & \dots        & a            & *
\end{array}
\right)\in \SL_n(A) | a\in A,\, a^{p^n}=0\}
$$
There is a natural
equivariant morphism
$$
\varphi:\P(V)\times \P(V^*) \rightarrow \P(V)\times \P(V^*)^{(r)}
$$
and $X$ is the scheme theoretic image $\varphi(Y)$. The induced
morphism $\varphi:Y\rightarrow X$ is the natural morphism
$$
G/P\rightarrow G/\tilde{P}
$$
given by the inclusion $P\subseteq \tilde{P}$.
Since the sheaf of ideals of the scheme theoretic image is the kernel
of the comorphism of $Y\rightarrow \P(V)\times \P(V^*)^{(r)}$
(\cite{Hartshorne},
Exercise II.3.11 (d)), $X$ is the zero scheme of the section $\bar{s}=
x_0^{p^r} \bar{y_0} +\dots + x_n^{p^r} \bar{y_n}$ of
$\O(p^r)\times\bar{\O}(1)$.
Using the isomorphism from \ref{fsubcover}, we get that $X=Z(\bar{s})$
is isomorhic to its scheme theoretic image $Z(\tilde s)\subseteq
\P(V)\times \P(V^*)$, where
$$
\tilde s=x_0^{p^r} y_0 + \dots x_n^{p^r} y_n
$$
is a section of $\O(p^r)\times \O(1)$.

\subsection{Cohomology of effective line bundles}

Let $a,b \in \Z$. The restriction to $Y$ of the line bundle $
\O(a)\times\O(b)$ on $\P^n\times \P^n$ will be denoted $L(a, b)$. The
restriction to $X$ of the line bundle $\O(a)\times \bar{\O}(b)$ on
$\P^n\times (\P^n)^{(r)}$ will be denoted $L(a, \bar{b})$.
Notice that the isomorphism in \ref{fsubcover} maps $L(a, \bar{b})$ to
$L(a, b)$.
By \ref{unsepdivisor} there is an
exact sequence
$$
0\rightarrow \O(-p^r)\times \O(-1) \rightarrow \O_{\P^n\times \P^n}
\rightarrow \O_X\rightarrow 0
$$
For the line bundle $L=L(a, \bar{b})$ on $X$, we therefore get
the exact sequence
$$
0\rightarrow \O(a-p^r)\times \O(b-1) \rightarrow
\O(a)\times \O(b) \rightarrow L \rightarrow 0
$$
Now assume that $a,b\geq 0$ ($L(a,\bar{b})$ is effective). Then
tracing through the long exact sequence and using the K\"unneth formula,
we get $\H^i(X, L)=0$, if $1\leq i< n-1$ along with the following
exact sequences:
$$
0\rightarrow \H^0(\P^n, \O(a-p^r))\otimes \H^0(\P^n, \O(b-1))\rightarrow
\H^0(\P^n, \O(a))\otimes\H^0(\P^n, \O(b))\rightarrow \H^0(X, L) \rightarrow 0
$$
and
$$
0\rightarrow \H^{n-1}(X, L)\rightarrow
\H^n(\P^n, \O(a-p^r))\otimes \H^0(\P^n, \O(b-1))\rightarrow 0
$$
By Serre duality one has
$$
\H^{n-1}(X, L)\cong \H^0(\P^n, \O(p^r-a-n-1))\otimes \H^0(\P^n, \O(b-1))
$$
so that the higher cohomology of $L$ vanishes if $a>p^r-n-1$.

By the adjunction formula (\cite{Hartshorne}, Proposition II.8.20) we get
$\omega_X\cong \O(p^r-n-1)\times \O(-n)$.
A line bundle $L=L(a,\bar{b})$ on $X$ is ample if and and only if $a,b>0$.
Kodaira type vanishing (vanishing higher cohomology
for $L\otimes \omega_X$, where
$L$ is ample) for $X$ amounts to the fact that $a+p^r-n-1>p^r-n-1$, when
$a>0$.

The unseparated incidence variety admits a lifting to a flat $\Z$-scheme.
There are projective homogeneous spaces
for $SL_4$, which do not admit a lifting to a flat $\Z$-scheme
(\cite{HaLa:varunsep}, \S 6).

\section{Structure of projective homogeneous spaces}

A projective homogeneous $G$-space $X$ is determined through its
stabilizer group scheme $G_x$ at some closed point $x\in X$. Notice that
since $X$ is projective, Borel's fixed point theorem implies that
$G_x$ contains a Borel subgroup $B$. We introduce some more
notation.
Let $T$ be a maximal torus of $G$ contained in the Borel subgroup $B$. Denote
by $R=R(G, T)$ the roots of $G$ w.r.t. $T$.
Let the roots $R(B, T)$ of
$B$ be the positive roots $R^+$ in $R$ and $S\subseteq R^+$ the
simple roots of $R$. Let $X(T)$ be the characters of $T$ and $Y(T)$
the one parameter subgroups. The usual pairing
$X(T)\times Y(T)\rightarrow \Z$ is denoted $\<\cdot, \cdot\>$. The coroot
in $Y(T)$ corresponding to $\alpha\in R$ is denoted $\alpha^\vee$.
The
root subgroup corresponding to $\alpha\in R$ is denoted by $U_{\alpha}$.
Let $\{X_{\alpha}\}_{\alpha\in R}$,
$\{H_{\alpha}\}_{\alpha\in S}$ be a Chevalley basis for $\Lie(G)$.
The monomials
$$
\prod_{\alpha\in R^+} X_{-\alpha}^{(n'(\alpha))}\,
\prod_{\alpha\in S} \binom{H_\alpha}{m(\alpha)}\,
\prod_{\alpha\in R^+} X_{\alpha}^{(n(\alpha))}
$$
where $n'(\alpha), m(\alpha), n(\alpha)\in\N$, form
a basis for the $k$-algebra of distributions $\Dist(G)$ (\cite{Jantzen},
II.1.12).
%Unlike the envelopping algebra in the complex case, $\Dist(G)$ is an
%infinitely generated $k$-algebra.
Recall that a subgroup scheme $H\subseteq G$ is uniquely determined by
its subalgebra $\Dist(H)\subseteq \Dist(G)$. A $k$-basis
for $\Dist((U_{\alpha})_n)$ is given by $1, X_{\alpha}, X_{\alpha}^{(2)},
\dots, X_{\alpha}^{(p^n-1)}$.

\subsection{Parabolic subgroup schemes}
\label{parsub}

Let $\tilde P$ be a subgroup scheme containing $B$. Since $P=\tilde P_{\red}$
is a parabolic subgroup (the nil-radical is a Hopf ideal), it follows
that $\tilde P$ is a connected group scheme. In particular we get
for $\alpha\in R^-$ that $\Dist(\tilde P)\cap\Dist(U_\alpha)=
\Dist((U_\alpha)_{n_\alpha})$ for a suitable
$n_\alpha$, where $0\leq n_\alpha \leq \infty$ with the convention
$(U_\alpha)_\infty=U_\alpha$. The subalgebra $\Dist(\tilde P)$ is
determined completely by $(n_\alpha)_{\alpha\in R^-}$.
With the assumptions given in \S \ref{charassumption} on $p=\kar\, k$
it follows by
(\cite{HaLa:varunsep}, Proposition 1.6) that $\Dist(\tilde P)$ is uniquely
determined by $(n_\alpha)_{\alpha\in-S}$.

One can construct $\tilde P\neq G$ as follows:
The maximal parabolic subgroup $P(\gamma)$
corresponding to a simple root $\gamma\in S$ is the
parabolic subgroup with roots generated by $S\setminus\{\gamma\}$.
The parabolic subgroup $P$ is the intersection
$P=P(\alpha_1)\cap\dots\cap P(\alpha_m)$ for certain
simple roots $S'=\{\alpha_1,\dots,\alpha_m\}\subseteq S$. It is easy
to see that
$\tilde P\subseteq G_n\,P(\alpha_i)$, for $n$ sufficiently big.
Let $n_i$ be the minimal $n$ with this property. Then
$$
\tilde P=G_{n_1}\, P(\alpha_1)\cap\dots\cap G_{n_m}\, P(\alpha_m)
$$
In the notation above $\Dist(\tilde P)$ is determined uniquely
by $n_{\alpha_1}=n_1,\dots, n_{\alpha_m}=n_m$ and $n_\alpha=\infty$ if
$\alpha\not\in S'$.

\subsection{The action of $G$ on $\P(L(\lambda))$}

Recall that the simple $G$-representations are parametrized by
dominant weights $X(T)_+$.
Let $L(\lambda)$ denote the simple $G$-representation associated
with $\lambda\in X(T)_+$. Then $L(\lambda)$ is generated by a $B$-stable
line of (highest) weight $\lambda$.

\begin{prop}
Let $S=\{\alpha_1,\dots, \alpha_l\}$ and
$\nu_p$ denote the $p$-adic valuation, such that $\nu_p(0)=\infty$.
Let $L(\lambda)$ be the simple representation of highest weight $\lambda\in
X(T)_+$ and $n_i=\nu_p(\<\lambda, \alpha_i^\vee\>)$.
Then
the stabilizer of the $B$-stable
line $x\in \P(L(\lambda))$ for the natural action of $G$ is
$$
G_x=G_{n_1} P(\alpha_1)\cap \dots\cap G_{n_l} P(\alpha_l)
$$
\end{prop}
\begin{pf}
Let $v$ be a generator for $x\in \P(L(\lambda))$. We compute the
algebra of distributions $\Dist(G_x)$. For the
induced action of $\Dist(G)$ on $L(\lambda)$ we have
$$
\Dist(G_x)=\{X\in \Dist(G) | X\, v=0\}
$$
By \ref{parsub} it suffices to show for a simple root $\alpha_i\in S$, that
$$
X_{-\alpha_i} \,v= X_{-\alpha_i}^{(2)}\, v
=\dots=X_{-\alpha_i}^{(p^{n_i}-1)}\,v=0
$$
and
$$
X_{-\alpha_i}^{(p^{n_i})}\, v\neq 0
$$
where the last condition is void in the case $n_i=\infty$.
Since $L(\lambda)$ is a highest weight module generated by $v$ and
$\alpha_i$ is a simple root, it suffices to prove that
$$
X_{\alpha_i}^{(n)} X_{-\alpha_i}^{(n)}\, v=0
$$
to conclude that $X_{-\alpha_i}^{(n)}\, v=0$.
In $\Dist(G)$ we have the following commutation formula for $\alpha\in R^+$:
$$
X_\alpha^{(m)} X_{-\alpha}^{(n)}=\sum_{j=0}^{\min(m,n)}
X_{-\alpha}^{(n-j)} \binom{H_\alpha-m-n+2j}{j} X_\alpha^{(m-j)}
$$
 From this formula it follows that
$$
X_{\alpha_i}^{(n)} X_{-\alpha_i}^{(n)}\, v= \binom{\<\lambda, \alpha_i^\vee\>}
{n}\, v
$$
When $n_i=\infty$ it follows that $X_{\alpha_i}^{(n)}\, v=0$ for $n>0$. Assume
now that $n_i<\infty$.
Since $\binom{m}{n}\equiv 0\, (\mod\, p)$ if $0<n<p^r$
and $\not\equiv 0\,(\mod\, p)$ if $n=p^r$, when $\nu_p(m)=r$, the
result follows.
\end{pf}

\label{gaction}
\label{Llambda}

\subsection{Exceptional parabolic subgroup schemes}

\label{ex}

The action of $G$ on $\P(L(\lambda))$ gives a lot of examples of
exceptional parabolic subgroup schemes in characteristic $2$  -
parabolic subgroup schemes which are not the intersection of
thickenings of the maximal parabolic subgroups as in \S
\ref{parsub}. In this section $k$ is assumed to be of
characteristic $2$. Recall that the simple module $L(\lambda)$
is a quotient of the Weyl module $V(\lambda)$. The Weyl module
$V(\lambda)$ is the base extension to $k$ of the minimal admissible
$\Z$-form in the simple representation of highest weight $\lambda$ for the
complex semisimple Lie algebra corresponding to $G$. Let $K(G)$ denote
the Grothendieck group of $G$. In the examples below, the decompositions
in $K(G)$ of Weyl modules were computed using Jantzen's sum
formula (\cite{Jantzen}, II.8). Example \ref{abuch} was discovered
using a computer program, developed by A.~Buch, for computations in modular
representation theory.

\begin{example}

Let $G$ be of type $B_2$ with positive roots $\alpha$, $\beta$,
$\alpha+\beta$ and $2\alpha+\beta$, where $\beta$ is the long simple root.
Let $\omega$ be the fundamental weight dual to $\beta$. In
$K(G)$ we have
$$
V(\omega)=L(\omega)+L(0)
$$
Let $v$ be a highest weight vector of $L(\omega)$.
To determine $G_x$, where $x=k\,v\in L(\omega)$,
we notice that $X_{-\alpha}\, v=0$, $X_{-\alpha-\beta}\, v=0$ (this
is because $0$ is not a weight of $L(\omega)$),
$X_{-\alpha-\beta}^{(2)}\, v\neq 0$, $X_{-2\alpha-\beta}\, v\neq 0$.
This means in the notation of \ref{parsub} that $\Dist(G_x)$ is given
by $n_{-\alpha}=0$,
$n_{-\beta}=\infty$, $n_{-\alpha-\beta}=1$,
$n_{-2\alpha-\beta}=0$.

\end{example}

\begin{example}

\label{abuch}

Let $G$ be of type $C_4$ with simple roots and fundamental dominant
weights numbered as below

\begin{picture}(150,40)(-125,10)

\put(50,20){\circle{3}}
\put(80,20){\circle{3}}
\put(110,20){\circle{3}}
\put(140,20){\circle{3}}

\put(47,10){$\alpha_1$}
\put(77,10){$\alpha_2$}
\put(107,10){$\alpha_3$}
\put(137,10){$\alpha_4$}

\put(47,30){$\omega_1$}
\put(77,30){$\omega_2$}
\put(107,30){$\omega_3$}
\put(137,30){$\omega_4$}

\put(51.5,19.5){\line(1,0){27}}
\put(81.5,19.5){\line(1,0){27}}
\put(111,21){\line(1,0){28}}
\put(111,18.3){\line(1,0){28}}

\put(123,17){$<$}

\end{picture}

\noindent
In $K(G)$ we have
$$
V(\omega_4)=L(\omega_4)+L(\omega_2)+L(0)
$$
and furthermore $\dim\, L(\omega_4)=16$, while $\dim\, V(\omega_4)=42$.
Let $v$ be a highest weight vector in $L(\omega_4)$ and $x=k\, v\in
\P(L(\omega_4))$. The stabilizer $G_x$ is given by $\Dist(G_x)$, which
is determined by the table
$$
\begin{array}{llll}
\alpha\in R^+ &    n_{-\alpha}            &
\alpha\in R^+ &    n_{-\alpha}\\
1000 & \infty      &  1100 & \infty     \\
1110 & \infty      &  0100 & \infty     \\
0110 & \infty      &  0010 & \infty     \\
0001 & 0      &  0011 & 1     \\
0111 & 1      &  1111 & 1     \\
0021 & 0      &  0121 & 1     \\
1121 & 1      &  0221 & 0     \\
1221 & 1      &  2221 & 0
\end{array}
$$
%$$
%\begin{array}{llcllc}
%\alpha\in R^+ & weight\ coord. & \<\omega, \alpha^\vee\> &
%\alpha\in R^+ & weight\ coord. & \<\omega, \alpha^\vee\> \\
%1000 & (0\ ,0\ ,0\ ,0)  & 0 & 1100 & (0\ ,0\ ,0\ ,0) & 0\\
%1110 & (0\ ,0\ ,0\ ,0)  & 0 & 0100 & (0\ ,0\ ,0\ ,0) & 0\\
%0110 & (0\ ,0\ ,0\ ,0)  & 0 & 0010 & (0\ ,0\ ,0\ ,0) & 0\\
%0001 & (0\ ,0\ ,-2,2) & 1 & 0011 & (0\ ,-1,0\ ,1)& 2\\
%0111 & (-1,1\ ,-1,1)& 2 & 1111 & (1\ ,0\ ,-1,1)& 2\\
%0021 & (0\ ,-2,2\ ,0) & 1 & 0121 & (-1,0\ ,1\ ,0)& 2\\
%1121 & (1\ ,-1,1\ ,0) & 2 & 0221 & (-2,2\ ,0\ ,0)& 1\\
%1221 & (0\ ,1\ ,0\ ,0)  & 2 & 2221 & (2\ ,0\ ,0\ ,0) & 1
%\end{array}
%$$
The orbit $X=G/G_x$ of $x=[v]$ has dimension $10$ and we get an
example of a variety lying on the boundary of Hartshorne's
conjecture \cite{LazVen} ($10=\frac{2}{3}15$)
$$
X\hookrightarrow \P(L(\omega_4))=\P^{15}
$$
I do not know whether $X\subseteq \P^{15}$ is a complete intersection.
One may check in accordance with Zak's result \cite{LazVen} on linear
normality that the restriction map
$$
\H^0(\P^{15}, \O(1))\rightarrow \H^0(X, \O_X(1))
$$
is surjective.

\end{example}

\section{Line bundles}

\label{linebundles}

In this section we classify the line bundles on projective homogeneous
$G$-spaces following \cite{Lau:Linebundles}. When $G$ is simply connected, all
line bundles are homogeneous induced by a character on $G_x$.

\subsection{Characters}

Let $X$ be a projective homogeneous $G$-space. Suppose that $G_x$ is the
stabilizer group scheme at a closed point $x\in X(k)$. Let $B$ be the Borel
subgroup contained in $G_x$.
The character lattice $X(B)=X(T)$ is
$$
\Z \omega_{\alpha_1}+\dots+\Z \omega_{\alpha_l}
$$
where $\omega_{\alpha}$ is the fundamental dominant weight
associated with the simple root $\alpha\in S$.
The restriction homomorphism $X(H)\rightarrow X(T)$ is injective for any
subgroup scheme $H\supseteq T$.
Recall that for a maximal parabolic subgroup $P(\alpha)$, we have
$X(P(\alpha))=\Z\, \omega_{\alpha}$.

\begin{lem}
\label{parchar}
Let $\alpha\in S$ be a simple root. Then
$$
X(G_n\, P(\alpha))=\Z p^n \omega_{\alpha}
$$
\end{lem}
\begin{pf}
It follows from (\cite{Jantzen}, II.3.15, Remarks 2),
that a character on the $n$-th
Frobenius kernel $G_n$ has to be trivial. Now the first isomorphism
theorem for groups gives
\begin{eqnarray*}
X(G_n P(\alpha))&=& X(G_n P(\alpha)/G_n) = X(P(\alpha)/G_n\cap P(\alpha))\\
                &=& X(P(\alpha))/P(\alpha)_n)=X(P(\alpha)^{(n)})=p^n
X(P(\alpha))\\
                &=& \Z p^n \omega_{\alpha}
\end{eqnarray*}
\end{pf}

Let $x_\alpha:\G_a\rightarrow G$ be the root homomorphism associated
with $\alpha\in R$. There is a homomorphism (\cite{Jantzen}, p.~176)
$$
\varphi_{\alpha}:\SL_2\rightarrow G
$$
such that
$$
\varphi_{\alpha}\left(\begin{array}{ll} 1 & a \\ 0 & 1 \end{array}\right)=
x_{\alpha}(a)\,\,\text{and}\,\,
\varphi_{\alpha}\left(\begin{array}{ll} 1 & 0 \\ a & 1 \end{array}\right)=
x_{-\alpha}(a)
$$
and
$$
\alpha^\vee(t)=\varphi_\alpha\left(
\begin{array}{ll}
t & 0 \\
0 & t^{-1}
\end{array}
\right),\, t\in \G_m
$$

We are now ready to prove

\begin{prop}
Let $G_x=G_{n_1}\,P(\alpha_1)\cap\dots\cap G_{n_m}\, P(\alpha_m)$, where
$\alpha_1,\dots,\alpha_m\in S$. Then
$$
X(G_x)=\Z p^{n_1} \omega_{\alpha_1} + \dots + \Z p^{n_m} \omega_{\alpha_m}
$$
\end{prop}
\begin{pf}
As $G_x\subseteq G_{n_i}\,P(\alpha_i),\ i=1, \dots,m$, we get by
lemma \ref{parchar}
$$
X(G_x)\supseteq \Z p^{n_1} \omega_{\alpha_1} +
\dots + \Z p^{n_m} \omega_{\alpha_m}
$$
Suppose on the other hand that $\lambda=a_1 \omega_{\alpha_1}+\dots
+a_m \omega_{\alpha_m}\in X(G_x)$. Then $\lambda\circ \varphi_{\alpha}$
is a character of $\tilde B\subseteq \SL_2(k)$, where for any
$k$-algebra $A$
$$
\tilde B (A)=\{
\left(
\begin{array}{ll}
*        & *    \\
a        & *
\end{array}
\right)\in \SL_2(A) | a\in A,\, a^{p^n}=0
\}
$$
and $n=\<\lambda, \alpha^\vee\>$.
Therefore we get
$p^{n_i} | \<\lambda, \alpha^\vee\>$, so that $p^{n_i} | a_i$.
\end{pf}

\subsection{Ample line bundles}

\label{ample}

Let $\tilde P$ be a parabolic subgroup scheme and $\chi\in X(\tilde P)$. The
total space of the line bundle $L_{\tilde P}(\chi)$ induced by $\chi$ is
$G\times^{\tilde P} \A^1=G\times \A^1/\tilde P$, where $\tilde P$ acts on
$G\times \A^1$ through $h.(g, a)=(g\, h, \chi(h^{-1}) a)$.
The natural morphism
$$
G\times^{\tilde P}\A^1\rightarrow G/\tilde P
$$
is equivariant, when $G$ acts on $G\times^{\tilde P} \A^1$ through
left multiplication. Since $G$ is simply connected ($\Pic G=0$) every
line bundle is induced by a character.

Let $P=\tilde P_{\red}$.
A line bundle $L$ on $G/P$, where
$P=P(\alpha_1)\cap \dots \cap P(\alpha_r)$ for simple roots
$\alpha_1,\dots, \alpha_r\in S$, is induced by a character
$$
\chi\in X(P)=\Z\omega_{\alpha_1}+\dots+\Z\omega_{\alpha_r}\subseteq X(B)=X(T)
$$
Then $L=L_P(\chi)$ is very ample on $G/P$ if and only if it is
ample on $G/P$ if and only if $\chi\in X(P)^{++}=\{\lambda\in X(P)|
\<\lambda, \alpha_1^\vee\>>0,\dots, \<\lambda, \alpha_r^\vee\>>0\}$.
It is easy to see that $f^* L_{\tilde P}(\chi)=L_P(\chi)$, where
$f:G/P\rightarrow G/\tilde P$ is the natural morphism.

\begin{thm}
Let $X=G/\tilde P$ be a projective homogeneous space such that
$$
\tilde P\cong G_{n_1}P(\alpha_1)\cap \dots \cap G_{n_r} P(\alpha_r)
$$
where $\alpha_1,\dots, \alpha_r\in S$ are simple roots and
$n_1,\dots, n_r$ integers $\geq 0$.
Let $\chi=a_1 p^{n_1} \omega_{\alpha_1} +\dots +
a_r p^{n_r} \omega_{\alpha_r}\in X(\tilde P)$.
Then $L_{\tilde P}(\chi)$ is very ample
on $X$ if and only if $L_P(\chi)=f^*L_{\tilde P}(\chi)$ is very
ample on $G/P$, where
$P=P(\alpha_1)\cap \dots \cap P(\alpha_r)$ and $f$ is the
natural morphism $f: G/P\rightarrow G/\tilde P$.
\end{thm}
\begin{pf}
Consider the natural diagram
$$
\begin{CD}
  G/P @>>>  G/P(\alpha_1)\times \dots \times G/P(\alpha_r)\\
  @VfVV                 @VVV\\
  G/\tilde P @>j>> (G/P(\alpha_1))^{(n_1)}\times \dots \times
(G/P(\alpha_r))^{(n_r)}
\end{CD}
$$
By \ref{fkernel} and \ref{product} it follows that $j$ is a closed
immersion. Since $(G/P(\alpha_i))^{(n_i)}=G/G_{n_i} P(\alpha_i)
\cong G/P(\alpha_i)$ it follows that ample line bundles on
$(G/P(\alpha_i))^{(n_i)}$ are very ample. Since the natural
morphism $G/P(\alpha_i)
\rightarrow G/G_{n_i} P(\alpha_i)$ is a finite surjective morphism
it follows (\cite{Hartshorne}, Exercise III.5.7 (d)) that
$L_{G_{n_i} P(\alpha_i)}(a_i\, p^{n_i} \omega_{\alpha_i})$ is
very ample if and only if $a_i>0$. By the Segre embedding we have
that
$$
L=L_{G_{n_1}P(\alpha_1)}(a_1\, p^{n_1} \omega_{\alpha_1})\times
\dots \times L_{G_{n_r}P(\alpha_r)}(a_r\, p^{n_r} \omega_{\alpha_r})
$$
is very ample. Now that $j$ is a closed immersion and
$j^* L=L_{\tilde P}(\chi)$ the result follows.
\end{pf}

\newpage
\bibliographystyle{amsplain}
\ifx\undefined\bysame
\newcommand{\bysame}{\leavevmode\hbox to3em{\hrulefill}\,}
\fi

\end{document}